\documentclass{birkjour}
\usepackage[all]{xy}
\usepackage{euscript}
\usepackage[hypertex]{hyperref}
\usepackage[all]{xy}
\usepackage{amssymb,mathrsfs,euscript}
\usepackage{mathrsfs,euscript,amssymb,amsmath,bbold}
\newtheorem{theorem}{Theorem}[section]

\theoremstyle{definition}
\newtheorem{definition}[theorem]{Definition}
\newtheorem{example}[theorem]{Example}

\makeatletter

\def\Tc{{\T^c}}
\def\nc{{non-commutative}}
\def\compT{{\widehat \T}}
\def\T{{\mathbb T}}
\def\bn{{\bar\psi}}
\def\D#1{{D^{[#1]}}}
\def\Symmetric{{\Sym}}
\def\compS{{\widehat \Sym}}
\def\Lin{{\it Lin}}
\def\sign#1{(-1)^{#1}}\def\sgn#1{(-1)^\varepsilon}\def\Sym{{\mathbb S}}

\def\Ainfty{{$A_\infty$}}\def\Tc{{\mathbb T}^c}
\def\Linfty{{$L_\infty$}}\def\Sc{{\mathbb S}^c}\def\calL{{\mathcal L}}
\def\Boj{B\"orjeson}\def\zn#1{{(-1)^{#1}}}
\def\bfk{{\mathbb k}}\def\id{1\!\!1}\def\ot{\otimes}
\def\otexp#1#2{{#1}^{\otimes #2}}\def\epi{\to}\def\nic{\relax}
\def\Rada#1#2#3{#1_{#2},\dots,#1_{#3}}\def\rada#1#2{{#1},\ldots,{#2}}

\begin{document}

%-------------------------------------------------------------------------
% editorial commands: to be inserted by the editorial office
%
%\firstpage{1} \volume{228} \Copyrightyear{2004} \DOI{003-0001}
%
%
%\seriesextra{Just an add-on}
%\seriesextraline{This is the Concrete Title of this Book\br H.E. R and S.T.C. W, Eds.}
%
% for journals:
%
%\firstpage{1}
%\issuenumber{1}
%\Volumeandyear{1 (2004)}
%\Copyrightyear{2004}
%\DOI{003-xxxx-y}
%\Signet
%\commby{inhouse}
%\submitted{March 14, 2003}
%\received{March 16, 2000}
%\revised{June 1, 2000}
%\accepted{July 22, 2000}
%
%
%
%---------------------------------------------------------------------------

\def\listtodoname{List of Todos}
\def\listoftodos{\@starttoc{tdo}\listtodoname}
\makeatother
\title[Higher braces]{Higher braces via formal (non)commutative geometry}

\author[M.~Markl]{Martin Markl}

\address{Mathematical Institute of the Academy, {\v Z}itn{\'a} 25,
         115 67 Prague 1, Czech Republic}

\address{MFF UK, 186 75 Sokolovsk\'a 83, Prague 8, Czech Republic}

\email{markl@math.cas.cz}

\thanks{The author was supported by the Eduard \v Cech
  Institute P201/12/G028  \hfill\break and  RVO: 67985840.}


%\typeout{zmenit datum}
%\typeout{zmenit styl}
\date{\today}

%%% ----------------------------------------------------------------------

\begin{abstract}
  We translate the main result of~\cite{markl:ab} to the language of
  formal geometry. In this new setting we prove directly that the Koszul
  resp.~\Boj\ braces are pullbacks of linear vector fields over the
  formal automorphism \hbox{$\varphi(a) = \exp(a) -1$} in the Koszul,
  resp.\ $\varphi(a) = a(1-a)^{-1}$ in the \Boj\ case. We then argue
  that both braces are versions of the {\em same\/} object, once materialized
  in the world of formal commutative geometry, once in the \nc~one.
\end{abstract}

%%% ----------------------------------------------------------------------
\maketitle
%%% ----------------------------------------------------------------------

\allowdisplaybreaks

\bibliographystyle{plain}

\tableofcontents

\section*{Introduction}

Let $A$ be a graded commutative associative algebra with a degree
$+1$ differential $\nabla: A \to A$ which {need} \underline{not\/}  be a
derivation. {\em Koszul braces\/} are linear degree $+1$ maps
$\Phi_n^\nabla : \otexp An \to A$, $n \geq 1$, defined
by the formulas
\begin{align*}
\Phi^\nabla_1(a) &= \nabla(a),
\\
\Phi^\nabla_2(a_1,a_2) &= \nabla(a_1a_2)
 - \nabla(a_1)a_2 - \sign{\nic a_1a_2\nic }\nabla(a_2)a_1,  
\\
\Phi^\nabla_3(a_1,a_2,a_3) &=\nabla(a_1a_2a_3)
-\nabla(a_1a_2)a_3 
- \sign{\nic a_1\nic (\nic a_2\nic  + \nic a_3\nic)}  \nabla(a_2a_3)a_1 
\\
&\hphantom{=} \hskip .5em-\sign{\nic a_3\nic (\nic a_1\nic  + \nic a_2\nic )}\nabla(a_3a_1)a_2  + \nabla(a_1)a_2a_3
\\
&\hphantom{=} \hskip .5em+ 
\sign{\nic a_1\nic (\nic a_2\nic  
+ \nic a_3\nic )}\nabla(a_2)a_3a_1+\sign{\nic a_3\nic (\nic a_1\nic  +
  \nic a_2\nic )} \nabla(a_3)a_1a_2,
\\
& \hskip .7em \vdots    
\\
\Phi^\nabla_n(\Rada a1n) &= 
\sum_{1 \leq i \leq n} \sign{n-i} \sum_\sigma \varepsilon(\sigma)\nabla(a_{\sigma(1)}
  \cdots a_{\sigma(i)})a_{\sigma(i+1)} \cdots a_{\sigma(n)},
\end{align*}
for $a,a_1,a_2,a_3, \ldots\in A$.  The sum in the last line runs
over all $(i,n-i)$-unshuffles $\sigma$, and $\varepsilon(\sigma) =
\varepsilon(\sigma;\Rada a1n)$ is the Koszul sign (recalled below). As proved for
instance in~\cite{bering-damgaard-alfaro:preprint}, these operations form an
\Linfty-algebra.

The above construction, attributed to Koszul~\cite{Koszul}, is
sometimes called the {\em Koszul hierarchy\/}, cf.\
also~\cite{akman:preprint97,akman-ionescu,bering-damgaard-alfaro:preprint,bering:CMP07,voronov03:higher}. It
is, among other things, used to define higher order derivations of
commutative associative algebras: the operator $\nabla : A \to A$ is
an order $r$ derivation, $r \geq 1$, if $\Phi^\nabla_{r+1}$
vanishes. Higher order derivations are crucial for the author's
approach to the BRST complex of the closed string field
theory~\cite[Section~4]{markl:la}; for this reason we believe that the
braces might be interesting also for physicists.

The assumption of commutativity of $A$ is crucial.  Although the
operations  $\Phi^\nabla_{n}$  make sense for
general $A$, they do not form any reasonable structure if $A$
is not commutative.

The noncommutative analog of the Koszul hierarchy was found in
April 2013 by Kay \Boj~\cite{Boj} who also proved that the result is
an \Ainfty-algebra; we recall his braces
$\{b^\nabla_n\}_{n \geq 1}$ in Subsection~\ref{noncko}. 
Amazingly, \Boj's braces are very different from the Koszul
ones. While $\Phi^\nabla_n$ consists of $2^n\!-\!1$
terms, $b^\nabla_n$ is for each $n \geq 3$ the sum of $4$ terms only!

In~\cite{markl:ab} we proved that both the Koszul and \Boj's braces
are the twistings of a linear \Linfty- (resp.~\Ainfty-) algebra
determined by $\nabla$. In this note we translate this result to the
language of formal (non)commutative geometry where \Linfty-
resp. \Ainfty-algebras appear as homological vector fields. Namely, we show 
in Theorems~\ref{Pakousove}
and~\ref{Pakousove-bis}  that both braces are
pullbacks of the linear vector field $\nabla$ 
over a formal automorphism $\varphi$  of the formal graded
affine pointed manifold $A$. This automorphism equals
\begin{subequations}
\begin{equation}
\label{eq:1}
\varphi(a) := \exp(a) -1 
\end{equation}
in the commutative (Koszul) case, and
\begin{equation}
\label{eq:2}
\varphi(a) := \frac a{1-a}
\end{equation}
\end{subequations}
in the \nc\ (\Boj) case. An important fact is that the Taylor
coefficients of both automorphisms are encapsulated in the
\underline{same\/} formal~sum
\[
\id_A + \mu^{[2]} +  \mu^{[3]} + \cdots  
\]
where $\mu^{[n]} : \otexp An \to A$, $n \geq 2$, is the iterated
multiplication in $A$. When interpreted in formal commutative
geometry, the associated formal automorphism is~(\ref{eq:1}) while
\nc\ geometry it is~(\ref{eq:2}).  The Koszul resp.~\Boj's braces are
thus versions of the \underline{same\/} object! 

As a pleasant by-product of our calculations, we get an important
general formula~(\ref{Musim_napsat_ten_uvodnik}) stated without proof
in~\cite{markl:ab}.

\hskip .5em

\noindent 
{\bf Conventions.}
All algebraic objects will be considered over a field $\bfk$ of
characteristic zero.  By $\id_U$ or simply by $\id$ when $U$ is
understood we denote the identity automorphism of a vector
space~$U$. We will reserve the symbol $\mu : A \ot A \to A$ for the
multiplication in an associative algebra $A$. The product $\mu(a,b)$
of elements $a,b \in A$ will usually be abbreviated as $ab$.

A {\em grading\/} will always mean a ${\mathbb Z}$-grading though 
most of our results  easily translate to the ${\mathbb
Z}_2$-graded setting, i.e.\ to the super world.
To avoid problems with dualization, we will  assume that all
graded vector spaces are of {\em finite type\/}.\footnote{A
  general case can be controlled by  a linear topology, but it would
  reach beyond the scope of the present paper.} For graded
vector spaces $U$ and $V$ we denote by $U\ot V$ their tensor product
over $\bfk$, and by $\Lin(U,V)$ the space of
degree $0$ $\bfk$-linear maps $U \to V$.

A permutation $\sigma\in\Sigma_n$ and graded variables $\Rada
\nu 1n$ determine the
\emph{Koszul sign} $\varepsilon(\sigma;\Rada   \nu 1n) 
\in \{-1,+1\}$ via the equation
\[
\nu_1\wedge\ldots\ \wedge
\nu_n=  \varepsilon(\sigma;\Rada   \nu1n) \cdot  \nu_{\sigma(1)}
\wedge\ldots\wedge
\nu_{\sigma(n)},
\]
in the free graded commutative associative
algebra ${\mathbb S}(\Rada \nu1n)$ generated by $\Rada \nu1n$. We usually
write  $\varepsilon(\sigma)$ instead of
$\varepsilon(\sigma;\Rada   \nu1n) $ when the meaning of $\Rada   \nu1n$
is clear.
For integers $u,v \geq 0$, an $(u,v)$-{\em unshuffle\/} is a 
permutation $\sigma \in \Sigma_{u+v}$ satisfying
\[
\mbox {
$\sigma(1)<\cdots<\sigma(u)$ \ and \ $\sigma(u+1)< \cdots < \sigma(u+v)$.}
\]

\section{Recollection of formal geometry}

In this section we recall basic concepts of formal (commutative)
geometry. We start with polynomial algebras, their completions and
duals, and explain how they are related to Taylor series of formal
maps. We then interpret \Linfty-algebras as homological
vector fields on formal graded pointed affine manifolds.
All notions recalled here are
standard and appeared in various forms in the literature. 
Their non-commutative variants are  briefly addressed in
Subsection~\ref{noncko}. 

\subsection{Algebras and coalgebras of symmetric tensors}

For a finite dimensional vector space $X$ denote by  
\[
\Sym(X) = \bigoplus_{k \geq 1} \Symmetric^k(X)
\]
the symmetric (polynomial) algebra generated by $X$. To distinguish
the multiplication in $\Sym(X)$ from other products that may occur in this
note we denote the product of two polynomials $p,q \in \Sym(X)$ by $p
\odot q$.

Let $A = X^*$ be the linear dual of $X$. One has, for each $n \geq 1$,
a non-degenerate pairing
\begin{equation}
\label{eq:25}
\langle - \, | \,  - \rangle : \Symmetric^n(X) \ot \Symmetric^n(A) \to \bfk
\end{equation} 
given by
\begin{equation}
\label{eq:22}
\langle x_1 \odot \cdots \odot x_n\, |\,  a_1 \odot \cdots \odot a_n
\rangle
:= \sum_{\sigma}  \ x_1(a_{\sigma(1)}) \cdots x_n(a_{\sigma(n)}),
\end{equation}
where $\Rada x1n \in X$, $\Rada a1n \in A$,  
and the summation ranges over all permutations $\sigma  \in \Sigma_n$.
As a particular case of~(\ref{eq:22}) we get
\begin{equation}
\label{eq:24}
\langle x_1 \odot \cdots \odot x_n \, | \, a \odot \cdots \odot a
\rangle
:=  n!\ x_1(a) \cdots x_n(a).
\end{equation}
Notice the factorial $n!$ emerging there.

We will need also the completion 
\[
\compS(X) = \prod_k \Sym^k(X)
\]
of $\Symmetric(X)$; $\compS(X)$ is the algebra of power series in $X$.  The
pairing~(\ref{eq:22}) extends to a non-degenerate pairing 
\begin{equation}
\label{eq:20}
\langle - \, | \, - \rangle : \compS(X) \ot
\Symmetric(A) \to \bfk
\end{equation}
 that identifies $\compS(X)$ with the linear dual of
$\Symmetric(A)$. This identification is such that the coalgebra structure on
$\Symmetric(A)$ induced from the algebra structure of $\compS(X)$
is given by the deconcatenation coproduct $\Delta : \Symmetric(A) \to
\Symmetric(A) \ot \Symmetric(A)$  
defined as
\begin{equation}
\label{eq:28}
\Delta(a_1\odot\cdots\odot a_n)=\sum_{1\leq j\leq n-1}\sum_\sigma
(a_{\sigma(1)}\odot\cdots\odot a_{\sigma(j)})\otimes
(a_{\sigma{(j+1)}}\odot\cdots\odot a_{\sigma(n)}),
\end{equation}
where $\Rada a1n \in A$ and $\sigma$ runs through all $(j,n-j)$ unshuffles.
We denote
$\Symmetric(A)$ viewed as a coalgebra with this coproduct by $\Sc(A)$.

The algebra  $\compS(X)$ is the free complete commutative associative
algebra generated by $X$. Therefore, for each complete commutative
associative algebra $B$ and a linear map $\omega : X \to B$, there
exists a unique morphism $h : \compS(X) \to B$ of complete
algebras such that the diagram
\[
\xymatrix@C = +4em@R = 1.2em{\compS(X) \ar@{-->}[r]^h 
& B 
\\
 \ar@/_/[ru]^\omega \rule{0em}{1em} X \ar@{_{(}->}[u]^\iota &
}
\]
in which $\iota : X  \hookrightarrow \compS(X)$
is the obvious inclusion, commutes.  In particular, an 
endomorphism $\phi : \compS (X) \to \compS (X)$ is  determined by the
composition
\[
X \stackrel\iota \hookrightarrow
\compS (X)  \stackrel \phi\longrightarrow \compS (X).
\]

Likewise, the coalgebra $\Sc(A)$ is the cofree
nilpotent cocommutative coassociative 
coalgebra cogenerated by~$A$ 
\cite[Definition~II.3.72]{markl-shnider-stasheff:book}. This means the
following.  Let \hbox{$\pi : \Sc (A) \twoheadrightarrow A$} be the obvious
projection. Then for any
nilpotent cocommutative coassociative 
coalgebra $C$ and for any linear map $\rho:C
\to A$, there exists exactly one coalgebra morphism $g:C \to
\Sc (A)$ making the diagram
\[
\xymatrix@C = +4em@R = 1.3em{C \ar@{-->}[r]^g \ar@/_/[rd]^\rho & \Sc (A) 
\ar@{>>}[d]^\pi
\\
&A
}
\]
commutative. In particular, any endomorphism $\phi : \Sc (A) \to
\Sc (A)$ is determined by the composition
\[
\xymatrix@C = 1.4em@1{\Sc (A) \ar[r]^\phi &\
\Sc (A) \ar@{>>}[r]^(.6)\pi &\ A
}.
\]

\subsection{Morphisms} 
\label{Jarka_opet_v_koncich}

Let $A$ be a finite dimensional 
vector space, $A^*$ its linear dual,  $\Sym(A^*)$ the symmetric 
algebra generated by
$A^*$, and $V$ another vector space. Elements
of $\Sym(A^*) \ot V$ are polynomials with coefficients in $V$. Each
 $p \in \Sym(A^*) \ot V$
determines a (non-linear) map $f: A \to V$ as follows.

The pairing~(\ref{eq:25}) in the obvious manner extends to a linear map
\[
\langle - \, | \, - \rangle : \Symmetric^n(A^*) \ot V \ot
\Symmetric^n(A) \longrightarrow V.
\]
We associate to every homogeneous polynomial $p_n \in \Sym^n(A^*) \ot V$, $n \geq 1$, a~map
$f_n : A \to V$ defined as
\begin{equation}
\label{Jarka_mne_zas_trapi}
f_n(a) : = \big\langle p_n \, | \, \D n(a) \big\rangle \in V ,\ a \in A,
\end{equation}
where $\D n : A \to \Symmetric^n(A)$ is the `diagonal' given by
\[
\D n(a) := \frac 1{n!} 
(\ \underbrace{a \odot  \cdots \odot a}_{n \ \hbox {\scriptsize
    times}}\ ) . 
\]
The factorial $n!$ in the definition of $\D n$ compensates the
factorial in~(\ref{eq:24}).
Every polynomial $p \in \Sym(A^*) \ot V$ is a finite sum 
\[
p = p_1+p_2+p_3+\cdots, \ p_n \in \Sym^n(A^*) \ot V,\ n \geq 1,
\]
of its homogeneous components.
We define $f : A \to V$ corresponding to $p$ as the finite sum
\begin{equation}
\label{eq:19}
f =f_1 + f_2 +f_3+ \cdots, 
\end{equation}
where $f_n$'s are, for $n \geq 1$, as in~(\ref{Jarka_mne_zas_trapi}).
If  the vector space $V$ equals the ground field $\bfk$,
$\Sym(A^*) \ot \bfk$ is isomorphic to $\Sym(A^*)$ and we obtain the standard
identification of the polynomial ring $\Sym(A^*)$ with the algebra of
regular functions on the affine space $A$.

We may formally extend the above construction to the space
\hbox{$\compS(A^*) \ot V$} of power series with coefficients in $V$. The sum
in (\ref{eq:19}) corresponding to \hbox{$p \in \compS(A^*) \ot V$}
then may be infinite.  In this way we interpret power series in
$\compS(A^*) \ot V$ as {\em Taylor coefficients\/} 
of maps $A \to V$.\footnote{Since
  we work formally, we do not pay any attention to the convergence
  issue.}

\begin{example}
\label{taylor}
Let $\bfk = {\mathbb R}$, $A$  be the real affine line  ${\mathbb R}$ with the
basic vector $e := 1\in  A$, and
 $x \in  A^*$ be such that $x(e) = 1$. For  a smooth map $\varphi
: A \to  A$ with $\varphi(0) = 0$ consider the power series
\[
p: =  \varphi'(0)\, x + \frac { \varphi''(0)}{2!}\, x^2 + \frac {\varphi'''(0)}{3!}\,  x^3 +
\cdots \
\in \compS(A^*)
\]
in $x$.  One then easily verifies that $p \in \Sc(A^*)$ determines,
via the above construction, a formal map $f: A \to A$ given by
\begin{equation}
\label{Jarka_by_zitra_mela_prijet_do_Mnichova}
f(\alpha e) =   \varphi'(0)\, \alpha + \frac { \varphi''(0)}{2!}\, \alpha^2 + \frac {\varphi'''(0)}{3!}\,  \alpha^3 +
\cdots \ , \ \alpha \in {\mathbb R},
\end{equation}
i.e.\ the Taylor expansion at $0$  of the smooth map $\varphi$.
\end{example}

An equivalent representation of formal maps $A \to V$ is based on the
identification
\begin{equation}
\label{eq:26}
\compS(A^*) \ot V \cong \Lin\big(\Sc(A), V\big)
\end{equation}
induced by the extended pairing~(\ref{eq:20}).
If $q : \Sc(A) \to V$ corresponds, under~(\ref{eq:26}), to $p \in \compS(A^*) \ot
V$   
then $f$ in~(\ref{eq:19}) is, for $n \geq 1$, simply the composition
\begin{equation}
\label{Jarka_mne_zas_trapi1}
f= q \circ(\id + \D 2 + \D 3 + \cdots).
\end{equation}

\begin{example} (Continuation of Example~\ref{taylor})
Let $A$ be again the real affine line ${\mathbb R}$ with the basic
vector $e = 1 \in {\mathbb R}$. Define $q : \compS(A) \to A$ by
\[
q(\ \underbrace{e \odot  \cdots \odot e}_{n \ \hbox {\scriptsize
    times}}\ ) := \varphi^{(n)}(0),
\]
where $\varphi^{(n)}(0)$ is, for $n \geq 1$, the $n$th derivative
of $\varphi$ at $0$. Clearly
\begin{equation}
\label{eq:3}
q = \varphi'(0)\, \id +  { \varphi''(0)}\, \mu^{[2]} +
{\varphi'''(0)}\, \mu^{[3]} + \cdots,
\end{equation} 
where \hbox{$\mu^{[n]}\! :\! \Symmetric^n(A) \to A$} denotes the
iterated multiplication in $A = {\mathbb R}$. With this
choice,~(\ref{Jarka_mne_zas_trapi1}) gives the Taylor
expansion~(\ref{Jarka_by_zitra_mela_prijet_do_Mnichova}). Notice
that~(\ref{eq:3}) does not contain any factorials.
\end{example}

\begin{example}
\label{sec:morphisms}
Suppose that $A$ is a commutative associative algebra and
\hbox{$V = A$}.  For each $n \geq 1$ one has the map 
\[
\mu^{[n]} :  \Sym^n(A) \to A
\]
that takes  $a_1 \odot \cdots \odot a_n \in  \Sym^n(A)$
to the iterated product $a_1 \cdots a_n \in A$.
Consider the morphism
\begin{equation}
\label{eq:29}
q  := \id + \mu^{[2]} +  \mu^{[3]} + \cdots:  \Sc(A) \to A.
\end{equation}
Summation~(\ref{Jarka_mne_zas_trapi1}) is then  the formal map
\begin{equation}
\label{eq:21}
\exp(a) -1 = a + \frac1{2!}a^2 +  \frac1{3!}a^3 + \cdots\ , \ a \in A, 
\end{equation}
whose inverse clearly equals
\begin{equation}
\label{eq:23}
 \ln(a+1) = 
\sum_{k \geq 1} \frac {(-1)^{k+1}}k\ a^k = 
a - \frac {a^2}2 + \frac {a^3}3 - \frac {a^4}4 + \cdots.
\end{equation}
\end{example}

\subsection{Formal manifolds}
\label{Jarusko-jarusko}

The previous constructions generalize to graded vector spaces. Then
$A$ is interpreted as a formal graded pointed affine
manifold,\footnote{We will sometimes omit `pointed' or/and `formal' to
  simplify the terminology.} with $\compS(A^*)$ playing the r\^ole of
its ring of functions.  The power series in $\compS(A^*) \ot V$
represent morphism from $A$ to $V$, where $V$ is interpreted as
another formal graded manifold.

Since we work with graded objects, some expressions should acquire
signs dictated by the Koszul sign convention. For instance, the terms
of the sum in the right hand side of~(\ref{eq:22}) must be multiplied
by the Koszul sign $\varepsilon(\rho)$ of the permutation
\[
\rho: 
\Rada x1n,\Rada a1n \longmapsto x_1, a_{\sigma(1)},\ldots,  
x_n,a_{\sigma(n)}.
\]
This is standard so we will pay no special attention to it.  Some
issues related to dualization may also arise, but they are not
the subject of this note.

Recall that a {\em vector field\/} $\chi$ on a (classical) smooth
manifold $M$ is a~section of its tangent bundle. The crucial property
of vector fields is that the assignment $f \mapsto \chi(f)$ is a
derivation of the ring of smooth functions $f: M \to {\mathbb R}$.
There is a contravariant action $\chi \mapsto \varphi^*\chi$ of
diffeomorphisms $\varphi : M \to M$ on vector fields characterized by
the formula
\begin{equation}
\label{eq:27}
\varphi^*\chi(f)\big(\varphi(a)\big)  = \chi(f \circ \varphi)(a)
\end{equation}
in which $f : M \to {\mathbb R}$ is a smooth function and $a \in M$ a
point.  We call $\varphi^*\chi$ the {\em pullback\/} of the vector
field $\chi$ {\em over\/} the diffeomorphism $\varphi$. 
Let us translate the pullbacks to the situation when we consider
instead of $M$ a formal affine graded manifold $A$
as follows.

Vector fields on $A$
are defined as 
derivations of the algebra of formal functions on $A$, i.e.\ derivations $\chi$ of
the algebra $\compS(A^*)$. 
An automorphism $\varphi:A \to A$ is given by 
a power series $p \in \compS (A^*) \ot A$, conveniently
represented using the duality  
by a map $\bar\phi : A^* \to \compS
(A^*)$. The invertibility of $\varphi$ is equivalent to the
invertibility  of the composition
\[
\xymatrix@1@C=1.5em{A^* \ar[r]^(.4){\bar\phi} 
& \ \compS (A^*) \
\ar@{>>}[r]^(.6)\pi
&\ A^*.
}
\]
Using the universal property of $\compS(A^*)$, we uniquely
extend $\bar\phi$ to an algebra automorphism $\phi : \compS(A^*) \to \compS(A^*)$.
If $\chi$ is a vector field on $A$, i.e.\ a derivation of
$\compS(A^*)$, the action~(\ref{eq:27}) is translated to the adjunction 
 \[
\varphi^*\chi :=  \phi^{-1} \circ \chi \circ \phi.
\]

\subsection{$L_\infty$-algebras as homological vector fields}
\label{sec:l_infty-algebras-as-1}

Recall that   an
\Linfty-algebra, also called sh or strongly
  homotopy Lie algebra, is a graded vector space  $L$ equipped with
  linear graded antisymmetric maps
$\lambda_n : \otexp Ln \to L$, $n \geq 1$, $\deg(\lambda_n) = 2-n$, 
that satisfy a set of axioms saying that $\lambda_1$ is a
differential, $\lambda_2$ obeys the Jacobi identity up to the homotopy
$\lambda_3$, etc., see for instance
\cite{lada-markl:CommAlg95} or~\cite{lada-stasheff:IJTP93}.

It will be convenient for the purposes  of this note to work with the version
transferred to the desuspension $A :=\ \downarrow\! L$. 
The transferred  $\lambda_k$'s have degree 
$+1$, are graded symmetric, and satisfy, 
for each $\Rada a1n
\in A$, $n \geq 1$, the axiom
\begin{equation}
\label{eq:8blis}
\sum_{i+j=n+1}\sum_\sigma
\varepsilon(\sigma)\lambda_j\big(\lambda_i(\Rada a{\sigma(1)}{\sigma(i)}),\Rada
a{\sigma(i+1)}{\sigma (n)}\big)=0,
\end{equation}
where $\sigma$ runs over $(i,n-i)$-unshuffles. In this guise,
\Linfty-algebra is an object $\calL = (A,\lambda_1,\lambda_2,\lambda_3,\ldots)$
formed by a graded vector space $A$ and degree $+1$ graded
symmetric linear maps $\lambda_n : \otexp An\to A$, $n \geq 1$,
satisfying~(\ref{eq:8blis}) for each $n \geq 1$.

Let $\Sc (A)$ be, as before, the symmetric coalgebra with the
deconcatenation coproduct~(\ref{eq:28}).  Thanks to the cofreeness of
$\Sc (A)$, each coderivation $\varrho$ of $\Sc (A)$ is 
determined by its components $\varrho_n := \pi \circ \varrho \circ
\iota_n$, $n \geq 1$, where $\pi: \Sc (A) \twoheadrightarrow A$ is the projection
and $\iota_n : {\mathbb S}^k( A) \hookrightarrow \Sc (A)$ the
inclusion. We write $\varrho =
(\varrho_1,\varrho_2,\varrho_3,\ldots)$.

In particular, let $\lambda := (\lambda_1,\lambda_2,\lambda_3,\ldots)$ 
be a degree $1$
coderivation of $\Sc (A)$ determined by the linear maps $\lambda_n$.  
By a classical result \cite[Theorem~2.3]{lada-markl:CommAlg95},
(\ref{eq:8blis}) shrinks to a single equation $\lambda^2 = 0$ .  Thus an
\Linfty-algebra is a degree $+1$ coderivation of $\Sc( A)$ that squares to zero.

The linear dual of $\lambda: \Sc(A) \to \Sc(A)$  is
a degree $-1$ 
derivation  $\vartheta$ 
of the algebra  $\compS(A^*)$ such that $\vartheta^2 = 0$. 
In the language of formal geometry,
$\vartheta$ is a degree $-1$ vector field on the formal affine
manifold $A$ that squares to zero. Such an object is called a
{\em homological vector field\/}. This is expressed by

\begin{definition}
\label{sec:l_infty-algebras-as}
\Linfty-algebras are homological vector fields on formal pointed graded affine
manifolds. 
\end{definition}

\begin{example}
\label{Vcera_mi_Jarunka_volala_kdyz_jsem_byl_v_Libni-bis}
Linear vector fields on the classical affine pointed manifold ${\mathbb R}^k$  
are those of the form
\[
f_1 \frac{\partial\hphantom{x^1}}{\partial x^1} + \cdots + f_k 
\frac{\partial\hphantom{x^k}}{\partial x^k},
\]
where $\Rada f1k :  {\mathbb R}^k   \to {\mathbb R}$ are linear
functions. The formal analog of linear vector fields are derivations
$\chi$ of $\compS(A^*)$ such that $\chi(A^*) \subset A^*$.

In the dual setting, linear vector fields are represented by
coderivations $\nabla$ of $\Sc(A)$ for which
the composition
\[
\xymatrix@1@C=1.5em{\Sym^{\geq 2}(A^*)\ \ar@_{^{(}->}[r] & \ar[r]^\chi
\ \Sc(A^*)\ & \ \Sc(A^*) \ \ar@{>>}[r]^(.65)\pi&\ A
}
\]
where $\Sym^{\geq 2}(A^*)$ is the subspace of $\Sc(A^*)$
of polynomials with vanishing linear term, is the zero map.
Such a coderivation is determined by its restriction to $A \subset
\Sc(A)$. This restriction is 
a linear map $A \to A$  which we denote $\nabla$ again.  
\end{example}

\begin{example}
\label{zase_krizova_epizoda_s_Jarkou}
An important particular case is  a linear homological vector field
given by a degree
$+1$ differential $\nabla : A \to A$, or
equivalently, by a `trivial' \Linfty-algebra  $(A,\nabla,0,0,\ldots)$.
\end{example}

\section{The braces}

We prove that the Koszul braces, as well as
their \nc\ analogs recalled in Subsection~\ref{noncko} below, are pullbacks
of linear vector fields over specific automorphisms. Throughout this section,
$A$ will be a graded associative, in
Subsection~\ref{dva_dny_s_Jarkou} also commutative, algebra, and
\hbox{$\nabla: A \to A$} a degree $+1$ differential which is not
required be a derivation. In Subsection~\ref{dva_dny_s_Jarkou} we
interpret $\nabla$ as a homological vector field on a (commutative)
formal manifold or, equivalently, as a trivial
\Linfty-algebra. In Subsection \ref{noncko} we interpret $\nabla$ as a
vector field on a \nc\ formal manifold, or as a trivial
\Ainfty-algebra, cf.~Examples~\ref{zase_krizova_epizoda_s_Jarkou}
resp.~\ref{zase_krizova_epizoda_s_Jarkou-bis} below.

\subsection{Koszul braces}
\label{dva_dny_s_Jarkou}
The aim of this subsection is to prove

\begin{theorem}
\label{Pakousove}
The Koszul braces recalled in the introduction are given by the
pullback of the linear homological field $\nabla$ over the formal
diffeomorphism $\exp - 1 : A \to A$.
\end{theorem}

We will study pullbacks of $\nabla$ over more general diffeomorphisms
and derive Theorem \ref{Pakousove} as a particular case. This 
generality would make the calculations more transparent. 
As a bonus we obtain formula~(\ref{Musim_napsat_ten_uvodnik}) given
without proof in~\cite{markl:ab}. 
Consider therefore a formal diffeomorphism
\begin{equation}
\label{mam_roztrhane_trenky}
A \ni  a 
\longmapsto \varphi(a) = f_1a + \frac{f_2}{2!}a^2 + \frac{f_3}{3!}
a^3 +\cdots,\ f_1,f_2,f_3 \in \bfk, \ f_1 \not= 0,
\end{equation}
having the inverse of the form
\[
A \ni  a  \longmapsto \psi(a) = g_1a + \frac{g_2}{2!}a^2 + \frac{g_3}{3!} a^3
+\cdots,  \ g_1,g_2,g_3,\ldots   \in \bfk.
\]
Notice that $g_1f_1 =1$.

As in Example~\ref{sec:morphisms}, one can easily check that
$\varphi$ is associated to
the map $q:  \Sc(A) \to A$ given by
\begin{equation}
\label{S_Jarkou_k_Pakousum}
q  := f_1 + f_2\mu^{[2]} +  f_3\mu^{[3]} + \cdots\ :  \Sc(A) \to A.
\end{equation}
Let us specify 
which map $\bar\phi : A^* \to \compS(A^*)$ corresponds to $q$ under
the extended pairing~(\ref{eq:20}). 

Since $q$ is an (infinite) sum of $f_n\mu^{[n]}$'s, we
determine first which map \hbox{$\bar\phi_n : A^* \to \Symmetric^n(A^*)$}
corresponds to $f_n\mu^{[n]} : \Symmetric^n(A) \to A$. Such a map
$\bar\phi_n$ is characterized by the duality
\[
\big\langle x \, | \,   f_n\mu^{[k]} (a_1 \odot \cdots
\odot a_n) \big\rangle = 
\langle \bar\phi_n(x) \, | \, a_1 \odot \cdots
\odot a_n \rangle
\] 
which has to hold for any $x \in A^*$ 
and $\Rada a1n \in A$. It is obvious that $\bar\phi_n$ must
equal $f_n\Delta^{[n]}$, with $\Delta^{[n]}$   
the iterated deconcatenation diagonal~(\ref{eq:28}),
\[
\Delta^{[n]} := (\Delta \ot \id^{\ot( n-2)})  \circ
(\Delta \ot \id^{\ot( n-3)}) \circ \cdots \circ \Delta
\]
where we put by definition $\Delta^{[1]} := \id_{A^*}$. We therefore have
by linearity
\[
\bar\phi = f_1  + f_2 \Delta + f_3 \Delta^{[3]} + 
f_4\Delta^{[4]} + \cdots.
\]
The formula for the map $\bar\psi: A^* \to \Sc(A^*)$ 
associated to the inverse
$\psi$ of $\varphi$ is analogous, namely
\[
\bar\psi = g_1  + g_2 \Delta + g_3 \Delta^{[3]} + 
g_4\Delta^{[4]} + \cdots.
\]
The map $\bar\phi : A^*
\to \compS(A^*)$ extends to a unique  automorphism
 $\phi$ of $\compS(A^*)$  given by
\begin{equation}
\label{eq:30}
\phi(x_1 \odot \cdots \odot x_n) =
\sum_{i_1,\ldots,i_n \geq 1} f_{i_1} \cdots f_{i_n} 
 \nabla^{[i_1]}(x_1) \odot
\cdots  \odot \nabla^{[i_n]}(x_n),
\end{equation}
$\Rada x1n \in A^*$.
There is a similar obvious formula for the extension of
$\bar\psi$, but we will not need it. 

It is however easier to work in the dual setting when vector fields
appear as {\em co}derivations of $\Sc (A)$ or, equivalently, as linear
maps $\Sc(A) \to A$. To see how $\varphi$ acts on vector fields in
this setup, we need to co-extend the map $q: \Sc(A) \to A$
in~(\ref{S_Jarkou_k_Pakousum}) to a coalgebra morphism $\Sc(A) \to
\Sc(A)$. We denote this co-extension by $\phi$ again, believing this
ambiguity will not confuse the reader. It acts on $a_1 \odot \cdots
\odot a_n\in \Symmetric^n(A)$ by\def\tau{\sigma}
\begin{eqnarray}
\lefteqn{
\label{Zase_krize_s_Jarkou}
\phi(a_1 \odot \cdots \odot a_n) =}
\\ \nonumber 
&=& \hskip -1em
\sum \epsilon(\tau) \frac{f_{i_1} \cdots f_{i_k}}{k!}  (a_{\tau(1)} \cdots
a_{\tau(i_1)})  
\odot \cdots \odot 
 (a_{\tau(i_1 + \cdots + i_{k-1} +1)} \cdots a_{\tau(n)}),
\end{eqnarray}
with the sum is taken over all $1 \leq k \leq n$, all $\rada
{i_1}{i_k} \geq 1$ such that $i_1 + \cdots + i_k = n$, and all
permutations $\tau \in \Sigma_n$ such that
\[
\tau(1) < \cdots < \tau(i_1),\ \ldots, 
\tau(i_1 + \cdots + i_{k-1} +1) < \cdots< \tau(n).
\]
Formula~(\ref{Zase_krize_s_Jarkou}) can be obtained either by
dualizing~(\ref{eq:30}) or directly, using the rule 
\[
\Delta\phi = (\phi \ot \phi) \Delta
\] 
describing the
interplay between morphisms of coalgebras and
coproducts.\footnote{Formula~(\ref{Zase_krize_s_Jarkou}) must be
  well-known, but we were unable to locate a  reference.} The
r\^ole of $k!$ in~(\ref{Zase_krize_s_Jarkou}) is explained in
Example~\ref{appear} below.

\vskip .5em
\noindent
{\bf Convention.}  To shorten the expressions, we will \underline{not}
write the Koszul signs explicitly as they can always be easily filled
in. \def\sqodot{{\hbox{$\,\odot\,$}}} We will also use the shorter
$\phi(a_1,\cdots, a_n)$ instead of $\phi(a_1 \odot \cdots \odot a_n)$,
etc.

\begin{example}
\label{appear}
Formula~(\ref{Zase_krize_s_Jarkou}) for $a,b,c \in A$ gives
\begin{align*}
\phi(a,\, & b, c)=
\\ 
= &f_3 (abc) + \frac{f_1f_2}{2!}(ab \sqodot c + bc
\sqodot a + ca \sqodot b + a \sqodot bc + b \sqodot ca + c \sqodot ab)
\\ 
&+ \frac{f_1^3}{3!} (a\sqodot b \sqodot c + b\sqodot c \sqodot a + c\sqodot a \sqodot b +
a\sqodot c \sqodot b + c\sqodot b \sqodot a + b\sqodot a \sqodot c)
\\
&= f_3 ( abc) + f_1f_2( ab \sqodot c + bc
\sqodot a + ca \sqodot b) + f_1^3  (a\sqodot b \sqodot c).
\end{align*}
The factorial in~(\ref{Zase_krize_s_Jarkou}) therefore removes
the multiplicities, so that each type of a term appears only once.
\end{example}

\begin{example}
\label{sec:koszul-braces}
Let us compute some initial terms of the composition $\bar\psi \nabla
\phi : \Sc(A) \to A$. For $a\in A$ we have
\[
\bar\psi\nabla \phi(a) = g_1f_1\nabla(a) = \nabla(a)
\]
where we used that $g_1f_1 =1$. For $a,b\in A$,
\begin{align*}
\bar\psi\nabla \phi(a, b)   &=
\bn \nabla (f_2 \ ab + f_1^2 a\sqodot b)
\\
& = \bn\big[f_2\nabla (ab) + f_1^2\big( \nabla
(a) \sqodot b +  a \sqodot \nabla(b)\big)\big]
\\
&= g_1f_2\nabla(ab) + g_2f_1^2\big (\nabla(a)b + a\nabla(b)\big),
\end{align*}
and, finally, for $a,b,c \in A$ one has
\begin{align*}
\bn \nabla \phi (a, b,c) & = 
\bn \nabla \big(f_3\,abc + f_2f_1( ab \sqodot c + bc
\sqodot a + ca \sqodot b) + f_1^3 \,  a\sqodot b \sqodot c\big)
\\
& =  \bn \Big[f_3\nabla(abc) +  f_2f_1 \big(\nabla(ab) \sqodot c +\nabla( bc)
\sqodot a + \nabla(ca) \sqodot b\big)
\\ 
&\hphantom{ =   \bn\varsigma^*(}+ f_2f_1\big( ab \sqodot \nabla(c) + bc
\sqodot \nabla(a) + ca \sqodot \nabla(b)\big)
\\ 
&\hphantom{ =   \phi{-1}\varsigma^*(}+f_1^3 \big (\nabla(a)\sqodot b \sqodot c + 
a\sqodot \nabla(b) \sqodot c + a\sqodot b \sqodot \nabla(c)\big)\Big]
\\
& =  g_1f_3 \nabla(abc) +g_2f_2f_1\big(\nabla(ab) c +\nabla( bc)a + 
\nabla(ca)b\big)
\\ 
&\hphantom{ = }\ +  (g_3f_1^3 + g_2f_2f_1)\big(\nabla(a) b c + 
a \nabla(b)  c +a  b\nabla(c)\big).
\end{align*}
\end{example}

\begin{example}
For $\varphi(a) = \exp(a)-1$, $\psi = \ln(a+1)$ one has
\[
f_1 = f_2 = f_3 = \cdots = 1
\]
and
\[ 
g_1 = 1,\ g_2 = -1,\ g_3 = 2,\ \ldots,\ g_n
= (-1)^{n-1}(n-1)!
\]
With this particular choice, the
calculations of Example \ref{sec:koszul-braces} lead, up to
implicit Koszul signs, to the Koszul braces $\Phi^\nabla_1$,
$\Phi^\nabla_2$  and $\Phi^\nabla_3$ recalled in the introduction,
i.e.
\begin{equation}
\label{Jarka_zitra}
\bar\psi\nabla \phi(a_1\odot \cdots \odot a_n) =
\Phi^\nabla_n(a_1\odot \cdots \odot a_n)
\end{equation}
for $n=1,2,3$.
\end{example}

Let us derive a general formula for $\bn \nabla \phi: \Sc(A) \to
A$. Using~(\ref{Zase_krize_s_Jarkou}),  one obtains
\begin{align*}
\bn \nabla& \phi(\Rada a1n) =
\\
&=
\bn \nabla \sum \frac{f_{i_1}\cdots f_{i_k}}{k!}  (a_{\tau(1)} \cdots
a_{\tau(i_1)})  
\odot \cdots \odot 
 (a_{\tau(i_1 + \cdots + i_{k-1} +1)} \cdots a_{\tau(n)})
\\
&=
\bn  \sum \frac{f_{i_1}\cdots f_{i_k}}{(k-1)!}  \nabla(a_{\tau(1)} \cdots
a_{\tau(i_1)})  
\odot \cdots \odot 
 (a_{\tau(i_1 + \cdots +i_{k-1} +1)} \cdots a_{\tau(n)})
\\
&=
  \sum  g_k\frac{f_{i_1}\cdots f_{i_k}}{(k-1)!}
 \nabla(a_{\tau(1)} \cdots
a_{\tau(i_1)}) 
 a_{\tau(i_{1} +1)} \cdots a_{\tau(n)}.
\end{align*}
The summations in the above display runs over the same data as
in~(\ref{Zase_krize_s_Jarkou}).
The substitution $i_1 \mapsto i$ converts the last line into
\def\doubless#1#2{{
\def\arraystretch{.5}
\begin{array}{c}
\mbox{$\scriptstyle #1$}
\\
\mbox{$\scriptstyle #2$}
\end{array}\def\arraystretch{1}
}}
\begin{equation}
\label{Musim_napsat_ten_uvodnik}
\bn \nabla \phi(a_1, \cdots,  a_n) = \hskip -.3em
\sum_{1 \leq i \leq n} \sum_\sigma c_{r} f_i
\nabla (\rada{a_{\sigma(1)}}{a_{\sigma(i)}})
a_{\sigma(i+1)} \cdots a_{\sigma(n)}
\end{equation}
where $\sigma$ runs over all
$(i,n-i)$-unshuffles, $r = n-i+1$, and
\begin{equation}
\label{eq:31}
c_r := \sum_{k \geq 2}
\ \sum_{\doubless{i_2 + \cdots +i_k = r}{\Rada i2k \geq 1}} 
g_k \frac{f_{i_2}\cdots f_{i_k}}{(k-1)!}
\frac{r!}{i_2! \cdots  i_k!} .
\end{equation}
The integer\def\kappa{\sigma}
\[
\frac{r!}{i_2! \cdots  i_k!} 
\] 
is, for $r = (i_2+ \cdots + i_k)$, 
the number of permutations $\sigma
\in \Sigma_n$ as in the sum~(\ref{Zase_krize_s_Jarkou}) with fixed
values $\kappa(1),\ldots,\kappa(i_1)$.

It is a simple exercise on manipulations with power series that
\[
c_{r} := \frac{d^r \psi'\big(\varphi(a)\big)}{da^r}\Big|_{t=0}.
\]
Formula~(\ref{Musim_napsat_ten_uvodnik}) is the one we gave  
without proof in~\cite[Section~2.4]{markl:ab}.

In the situation of Theorem~\ref{Pakousove},  
$\varphi(a) = \exp(a)-1$, $\psi(a) =
\ln(a+1)$, so  
$\psi'(a) = (1+a)^{-1}$, thus $\psi'\big(\varphi(a)\big) = e^{-a}$,
therefore $c_{r} = (-1)^r$ for $r \geq 1$. Since $f_s =1$ for each $s \geq 1$,
formula~(\ref{Musim_napsat_ten_uvodnik}) 
reproduces the Koszul braces in the introduction as direct inspection
shows, i.e.~(\ref{Jarka_zitra}) holds for every $n \geq 1$.
A combinatorial by-product of our calculations is 
the equality
\[
 \sum_{k \geq 2}
\ \sum_{\doubless{i_2 + \cdots +i_k = r}{\Rada i2k \geq 1}}
 (-1)^{k-1}
\frac{r!}{i_2! \cdots  i_k!}  = (-1)^r.
\]
We do not know any elementary proof of this fact.

\subsection{\Boj\ braces}
\label{noncko}

A non-commutative analog of Koszul braces was constructed 
by K.~\Boj. For a graded associative, not necessarily
commutative, algebra $A$ and a degree $+1$ differential \hbox{$\nabla:
  A \to A$}, he defined in~\cite{Boj} linear \hbox{degree $+1$} operators
\hbox{$b^\nabla_n : \otexp An \to A$}  by
\begin{align*}
b^\nabla_1(a) &= \nabla(a),
\\
b^\nabla_2(a_1,a_2) &= \nabla(a_1a_2) - \nabla(a_1)a_2 - \zn {a_1} 
a_1\nabla(a_2),  
\\
b^\nabla_3(a_1,a_2,a_3) &=\nabla(a_1a_2a_3) -\nabla(a_1a_2)a_3
\\
&  \hphantom{= }\ -
\zn {a_1}  a_1\nabla(a_2a_3)
+\zn {a_1}  a_1\nabla(a_2)a_3,
\\
b^\nabla_4(a_1,a_2,a_3,a_4) &=\nabla(a_1a_2a_3a_4)\!
 - \! \nabla(a_1a_2a_3)a_4
\\
&  \hphantom{= }\
\!-\!\zn {a_1}  a_1\nabla(a_2a_3a_4)\! +\!\zn {a_1}  a_1\nabla(a_2a_3)a_4,
\\
&\hskip .6em \vdots
\\
b^\nabla_k(\Rada a1k) &= 
\nabla(a_1\cdots a_k) - \nabla(a_1\cdots a_{k-1})a_k 
\\
&\hphantom{=} \hskip .5em-
\zn{a_1}a_1\nabla(a_2\cdots a_{k}) + \zn{a_1}a_1 \nabla(a_2\cdots a_{k-1})a_k.
\end{align*}
for $a,a_1,a_2,a_3,\ldots \in
A$.
He also proved that these operators form an
\Ainfty-algebra. In~\cite{markl:ab} we showed that \Boj\ braces are,
as their commutative counterparts, a
twisting of  $\nabla$ interpreted as a 
trivial \Ainfty-algebra. In this subsection we put these results into
the context of non-commutative \hbox{geometry}.

Everything in fact translates verbatim with only minor modifications
from the commutative case analyzed in the previous parts of this
note. For a finite-dimensional vector space $A$ we denote by
\[
\T(A^*) = \bigoplus_{k\geq 1} \T^k(A^*), \  
\T^k(A^*) := \underbrace{A^* \ot \cdots \ot A^*}_{k \ \hbox {\scriptsize
    times}}
\]
the tensor algebra of its dual. 
If $V$ is another finite-dimensional vector space, then every
homogeneous \nc\ polynomial \hbox{$p_n \in \T^n(A^*) \ot V$} for \hbox{$n \geq
1$} determines a~map $f_n : A \to V$ defined as
in~(\ref{Jarka_mne_zas_trapi}), with the only difference that 
the `diagonal' 
\[
\D n(a) :=
(\ \underbrace{a \odot  \cdots \odot a}_{n \ \hbox {\scriptsize
    times}}\ ) , 
\]
now does \underline{not} involve the factorial. Since every $p \in
\T(A^*) \ot V$ is a finite sum of its homogeneous components, we can
linearly extend the above construction and assign to each $p$ a map
$f: A \to V$.
Passing to the completion
\[
\compT(A^*) := \prod_{k\geq 1} \T^k(A^*),
\]
we interpret \nc\ power series in $\compT(A^*) \ot V$ as \nc\ Taylor
coefficients of maps $A \to V$. Since we have, by duality, an
isomorphism
\[
\compT(A^*) \ot V \cong \Lin\big(\Tc(A),V\big)
\]
where $\Tc(A)$ is the tensor coalgebra with the deconcatenation
coproduct, we may equivalently represent \nc\ Taylor coefficients of maps
\hbox{$f: A
\to V$} by linear morphisms $q: \Tc(A) \to V$. Let us give
a non-commuta\-tive analog of Example~\ref{sec:morphisms}.

\begin{example}
Let $A$ be an associative algebra. Then
the morphism
\begin{equation}
\label{eq:29bis}
q  := \id + \mu^{[2]} +  \mu^{[3]} + \cdots:  \Tc(A) \to A
\end{equation}
represents the \nc\ Taylor series
\begin{equation}
\label{V_pondeli_letim_do_Mnichova}
\frac a {1-a} = a + a^2 + a^3 + \cdots,\  a\in A
\end{equation}
whose inverse equals
\[
\frac a {1+a} = a - a^2 + a^3 - \cdots,\  a\in A.
\]
\end{example}

Formal \nc\ differential geometry is build upon the above classical
non-graded finite-dimensional affine spaces analogously as explained in
Subsection~\ref{Jarusko-jarusko} for the commutative case. That is,
$A$ is now a
graded vector space interpreted as a \nc\ formal affine pointed
manifold with $\T(A^*)$ its \nc\ ring of regular functions.
Given another formal \nc\ affine  manifold $V$, the
\nc\ power series in $\compT(A^*) \ot V$ represent formal maps $f: A
\to V$. By duality, the same formal maps can equivalently be
represented by linear morphisms $q : \Tc(A) \to A$. 

Vector fields on a formal \nc\ manifold $A$ are derivations of the complete
algebra $\compT(A^*)$ or, equivalently, coderivations of the tensor
coalgebra $\Tc(A)$. As in the commutative case, automorphisms act on
vector fields by adjunction.

\Ainfty-algebras are \nc\ versions of \Linfty-algebras recalled in
Subsection~\ref{sec:l_infty-algebras-as-1}, and their historical
precursors~\cite{stasheff:TAMS63}.  They consist of a 
graded vector space $U$
together with degree $2-n$ 
linear maps $m_n : \otexp Un \to U$, $n \geq 1$, such that
$m_1$ is a differential, $m_2$ is associative up to the homotopy
$m_3$, etc.

As in the case of \Linfty-algebras, we  transfer the
maps $m_n : \otexp Un \to U$ to the desuspension $A :=\
\downarrow\! U$. These transferred $m_n$'s are all 
of degree $+1$ and satisfy
\begin{equation}
\label{eq:8bis}
\sum_{u+v = n+1}\sum_{1\leq i \leq u}
m_u(\id_A^{\ot {i-1}} \ot m_v \ot \id_A^{\ot {u-i}}) = 0
\end{equation}
for each $n \geq 1$.

Since the tensor coalgebra
$\Tc (A)$ is the cofree nilpotent coassociative coalgebra
cogenerated by $A$, each coderivation $\varrho$
of $\Tc (A)$ 
is uniquely determined 
by its components $\varrho_n : \otexp An \to A$, $n \geq 1$, defined as
$\varrho_n := \pi \circ
\varrho \circ \iota_n$, where $\pi: \Tc (A) \epi A$ is the
projection and $\iota_n : \otexp An \hookrightarrow \Tc (A)$ the obvious
inclusion. We express this situation by  $\varrho =
(\varrho_1,\varrho_2,\varrho_3,\ldots)$. 

One in particular has a  degree $+1$ coderivation $m =
(m_1,m_2,m_3,\ldots)$   of $\Tc(A)$ determined by the \Ainfty-algebra
above. As in the \Linfty-case, the system~(\ref{eq:8bis}) is equivalent
to a single equation $m^2 = 0$. In other words, an \Ainfty-algebra is a
degree $+1$ coderivation of the tensor coalgebra $\Tc(A^*)$ that
squares to zero. Its linear dual $\varpi : \compT(A) \to \compT(A)$ is
a degree $-1$ derivation that squares to zero, i.e., a {\em homological
vector field\/}. 
We may therefore give the following analog of 
Definition \ref{sec:l_infty-algebras-as}:

\begin{definition}
  \Ainfty-algebras are homological vector fields on formal \nc\ graded
  pointed affine manifolds.
\end{definition}

The notion of linear vector fields translate verbatim from the
commutative case. Here is an analog of
Example~\ref{zase_krizova_epizoda_s_Jarkou}: 

\begin{example}
\label{zase_krizova_epizoda_s_Jarkou-bis}
A degree $+1$ differential $\nabla : A \to A$ extends to a degree $+1$
derivation of $\Tc(A)$, i.e.\ to a linear homological vector field on
the formal \nc\ affine pointed manifold $A$. In other words, it
determines  a~`trivial' \Ainfty-algebra $(A,\nabla,0,0,\ldots)$.
\end{example}

We finally formulate a \nc\ analog of Theorem~\ref{Pakousove}:

\begin{theorem}
\label{Pakousove-bis}
The \Boj\ braces  are given by the
pullback of the linear homological field $\nabla$ over the formal
diffeomorphism 
\[
\frac a{1-a} : A \to A
\] 
of the formal \nc\ pointed affine manifold $A$.
\end{theorem}

The proof of this theorem is simpler than that of  
Theorem~\ref{Pakousove}, since no symmetry enters. 
We leave it to the reader, as well as the \nc\ analog
of~(\ref{eq:31}) as appeared in~\cite[Section~2.2]{markl:ab}.

\def\cprime{$'$} \def\cprime{$'$}

%\bibliography{b}

\end{document}